\documentclass[letterpaper, 10 pt, conference]{ieeeconf}

\usepackage{amsthm}
\usepackage{amsmath,amssymb,amsfonts}
\usepackage{graphicx}
\usepackage{graphics}
\usepackage{varioref}
\usepackage{hyperref}
\usepackage[capitalise]{cleveref}
\usepackage{subfigure}
\usepackage{float}
\usepackage{cite}
\usepackage{algorithm}
\usepackage{algorithmic}
\usepackage{tikz}
\usepackage{comment}

\usepackage{soul}

\newcommand{\kl}{}
\newcommand{\ml}{}
\newcommand{\icl}{}
\newcommand{\il}{}
\newcommand{\ill}{}
\newcommand{\nkl}{}
\newcommand{\lmm}{}
\newcommand{\lk}{}

\usepackage[prependcaption,colorinlistoftodos]{todonotes}

\newcommand{\todoiny}[1]{}
\newcommand{\todoing}[1]{}
\newcommand{\todoKL}[1]{}
\newcommand{\todoinc}[1]{}

\IEEEoverridecommandlockouts                              
\title{\LARGE \bf Convergence Rate Bounds for the Mirror Descent Method:\\ IQCs and the Bregman Divergence}

\author{Mengmou~Li, Khaled~Laib, and Ioannis~Lestas
\thanks{\lmm{This work was supported by ERC starting grant 679774.}}
\thanks{The authors are with the Department of Engineering, University of Cambridge, Trumpington Street, Cambridge CB2 1PZ, United Kingdom.  Emails: \{ml995, kl507, icl20\}@cam.ac.uk}
}

\newtheorem{theorem}{Theorem}
\newtheorem{lemma}{Lemma}

\newtheorem{remark}{Remark}

\begin{document}

\maketitle
\thispagestyle{empty}
\pagestyle{empty}

\begin{abstract}
This paper is concerned with convergence analysis \ill{for the mirror descent (MD) method, a well-known algorithm in convex optimization.}
An analysis framework via integral quadratic constraints (IQCs) is {constructed} 
 to analyze the convergence rate of the MD method with strongly convex objective functions in both continuous-time and discrete-time.
We formulate the problem of finding convergence rates of the MD algorithms into feasibility problems of linear matrix inequalities (LMIs) in both schemes.
\icl{In particular, in continuous-time}, we show that \icl{the Bregman divergence function, which is commonly used as a Lyapunov function for this algorithm,} is a special case of \icl{the class of Lyapunov functions} \icl{associated with the} Popov criterion, \icl{when the latter is applied to an appropriate \ill{reformulation} of the problem}. Thus, applying the Popov criterion \icl{and its combination with other IQCs,}
\ill{can lead to convergence rate bounds with reduced conservatism.} 
\ml{We also \il{illustrate via examples that the 
convergence rate bounds derived \ill{can be} tight.}}
\todoing{\st{Last sentence not very clear}}
\end{abstract}
\todoing{\st{I think add all derivations in an appendix to improve the readability.}}

\section{Introduction}
The mirror descent (MD) method was initially proposed by Nemirovsky and Yudin \cite{nemirovskij1983problem} for solving constrained convex optimization problems. By choosing a Bregman distance function in place of the Euclidean distance to reflect the geometry of the constraint sets, it generalizes the gradient descent (GD) method from the Euclidean space to Hilbert and Banach spaces \cite{bubeck2014convex}.
Due to its applications in machine learning and large-scale optimization problems, it has received considerable research attention in various contexts, such as stochastic optimization \cite{duchi2012ergodic,nedic2014stochastic}, distributed optimization \cite{doan2018convergence,sun2021centralized}, and accelerated algorithms \cite{krichene2015accelerated,wibisono2016variational}.
\todoing{need to clarify for what classed of problems it is used\\
Sorry, I didn't get it. What do you mean by classes of problems?}

Many optimization algorithms can be treated as nonlinear dynamical systems, whose convergence may be verified by the Lyapunov stability theorem.
\ml{The 
Lyapunov function \ill{commonly used} for the MD method is the Bregman divergence function \lmm{measuring} the Bregman distance between the decision variable and the optimal solution.}
The Bregman divergence function was introduced by Bregman to find the intersection of convex sets \cite{bregman1967relaxation}. It has wide applications in the analysis of distributed optimization \cite{li2020input}, port-Hamiltonian systems \cite{jayawardhana2007passivity}, \ill{equilibrium independent stability\cite{simpson2018hill}, power systems \cite{de2017bregman,monshizadeh2019secant}, in addition to the MD method.}

\ill{Nevertheless, when bounds on the convergence rate need to be established it is important to have systematic methods that allow to construct Lyapunov functions with more advanced structures, or allow via other means to deduce convergence rates with reduced conservatism. It has been pointed out in the optimization literature that IQCs \cite{megretski1997system} can be a useful tool in this direction 
\cite{lessard2016analysis,dhingra2018proximal}.
However, their application in the case of the MD method is non-trivial as the MD dynamics involve the composition of two nonlinearities that correspond to monotone operators, with this composition not preserving these monotonicity properties.} 

\ill{
Our contributions in this paper can be summarized as follows:
\begin{enumerate}
	\item We show in continuous-time that the use of the Bregman divergence as a Lyapunov function for the MD method is a special case of Lyapunov functions that follow from the Popov criterion, when this is applied to an appropriate reformulation of the problem.
\item We use conic combinations of Popov IQCs and other type of IQCs that are relevant in our reformulation to derive convergence rate bounds for the MD method with reduced conservatism.
\end{enumerate}
}
\ill{
The convergence rate bounds deduced are formulated as solutions to LMIs in both discrete and continuous time. In the case of discrete time dynamics we also show via numerical examples that these bounds can be tight.}

The rest of this paper is organized as follows.
In \Cref{Preliminaries}, \icl{preliminaries on} the MD method and \icl{IQCs are provided.} The continuous-time and discrete-time MD methods are \icl{analysed} via IQCs in \Cref{Continuous-time mirror descent method} and \Cref{Discrete-time mirror descent method}, respectively.
In \Cref{Numerical Examples}, numerical examples are given to verify our results.
Finally, the paper is concluded in \Cref{Conclusion}.

\section{Preliminaries}\label{Preliminaries}
\subsection{Notation}
Let $\mathbb{R}$, $\mathbb{Z}$, $\mathbb{Z}_{+}$ denote the set of real numbers, integers, and nonnegative integers, respectively. Let $I_d$ and $0_d$ denote the $d \times d$ identity matrix and zero matrix, respectively. Their subscripts can be omitted if it is clear from the context. \ml{$\textup{diag} (\alpha_1, \ldots, \alpha_d)$ denotes a $d \times d$ diagonal matrix with $\alpha_i$ on its $i$-th diagonal entry.}
Let $\mathbf{RH}_{\infty}$ be the set of proper real rational functions without poles in the closed right-half plane. The set of $m \times n$ matrices with elements in $\mathbf{RH}_{\infty}$ is denoted $\mathbf{RH}_{\infty}^{m \times n}$. Let $\mathbf{L}_{2}^{m} [0, \infty)$ be the Hilbert space of all square integrable and Lebesgue measurable functions $f: [0, \infty) \rightarrow \mathbb{R}^{m}$. It is a subspace of $\mathbf{L}_{2e}^{m}[0, \infty)$ whose elements only need to be integrable on finite intervals. Let ${l}_{2}^{m}(\mathbb{Z}_{+})$ be the set of all square summable sequences $f : \mathbb{Z}_{+} \rightarrow \mathbb{R}^{m}$.
\ml{Given a Hermitian matrix $H (j\omega)$, $H^*(j\omega) : = H^T(-j\omega)$ represents its conjugate transpose and $\text{Re}\{ H (j\omega)\}$ denotes its real part.}
\todoKL{\st{Notation $\text{diag}$ is not defined. Also, use the same notation $\text{diag}$  or ${diag}$ (with or without the text environment).} }

Given $0 \leq \mu \leq L$, we denote $S(\mu, L)$ as the set of functions $f: \mathbb{R}^{d} \rightarrow \mathbb{R}$ that are continuously differentiable, $\mu$-strongly convex and $L$-smooth, i.e., $\forall x , ~ y$,
 \begin{align*}
 	\mu \| x - y \|^2 \leq \left( \nabla f (x) - \nabla f(y) \right)^T(x - y) \leq L \| x - y \|^2.
 \end{align*}
In this work, we assume $\mu > 0$ for all the functions we study if not specified otherwise. The condition number $\kappa$ of functions in $S(\mu, L)$ is defined by $\kappa := L/\mu \geq 1$.


\subsection{Integral quadratic constraints}
In continuous-time, a bounded operator $\Delta : \mathbf{L}_2^{n} [0, \infty) \rightarrow \mathbf{L}_{2}^{m} [0, \infty)$ is said to satisfy the IQC defined by $\Pi$, denoted by $\Delta \in \text{IQC}(\Pi)$, if
\begin{align}\label{eq:IQC definition}
	\int_{-\infty}^{\infty} \begin{bmatrix} \hat{v} (j\omega) \\ \hat{w} (j \omega) \end{bmatrix}^*
	\Pi(j\omega)
	\begin{bmatrix} \hat{v} (j\omega) \\ \hat{w} (j \omega) \end{bmatrix} d\omega \geq 0
\end{align}
for all $v \in \mathbf{L}_{2}^{n} [0, \infty)$ and $w = \Delta(v)$, where $\hat{v}(j\omega)$, $\hat{w} (j\omega)$ are the Fourier transforms of $v$, $w$, respectively, and $\Pi (j\omega)$ can be any measurable Hermitian valued function.
In discrete-time, condition \eqref{eq:IQC definition} is reduced to
\begin{align*}
	\int_{- \pi}^{\pi} \begin{bmatrix} \hat{v} ( e^{j\omega} ) \\ \hat{w} (e^{j \omega}) \end{bmatrix}^*
	\Pi (e^{j \omega})
	\begin{bmatrix} \hat{v} (e^{j\omega} ) \\ \hat{w} ( e^{j \omega} ) \end{bmatrix}
	d \omega \geq 0
\end{align*}
for all $ v \in {l}_{2}^{n}(\mathbb{Z}_{+})$, and $w = \Delta (v)$.

Define the truncation operator $P_T$ which does not change a function on the interval $[0, T]$ and gives the value zero on $(T, \infty]$.
The operator $\Delta$ is said to be \textit{causal} if $P_T \Delta P_T = P_T \Delta$, for all $T \geq 0$.
Consider the interconnection
\begin{equation}\label{eq:feedback interconnection model}
\begin{aligned}
	v = & Gw + g\\
	w = & \Delta (v) + e
\end{aligned}
\end{equation}
where $g \in \mathbf{L}_{2e}^{l}[0, \infty)$, $e \in \mathbf{L}_{2e}^{m}[0, \infty)$, $G$ and $\Delta$ are two causal operators on $\mathbf{L}_{2e}^{m}[0, \infty)$, $\mathbf{L}_{2e}^{l}[0, \infty)$, respectively. The feedback interconnection of $G$ and $\Delta$ is \textit{well-posed} if the map $(v, w) \mapsto (e, g)$ defined by \eqref{eq:feedback interconnection model} has a causal inverse on $\mathbf{L}_{2e}^{m+l}[0, \infty)$. The interconnection is \textit{stable} if, in addition, the inverse is bounded, i.e., there exists a constant $c > 0$ such that $\int_{0}^{T} \left( |v|^2  + |w|^2\right) d t \leq c \int_{0}^{T} \left( |g|^2 + |e|^2 \right) dt$.
System \eqref{eq:feedback interconnection model} with linear $G$ and \ml{static} nonlinear $\Delta$ is called \lmm{a} Lur'e system.

\todoing{\st{I think Lure system is about static nonlinearities}}
We will adopt the following IQC theorem for stability analysis.
\begin{theorem}[\hspace{1sp}\cite{megretski1997system}]\label{thm:IQC theorem}
	Let $G(s) \in \mathbf{RH}_{\infty}^{l \times m}$, and let $\Delta$ be a bounded causal operator. Assume that:
	\begin{enumerate}
		\item for every $\tau \in [0,1]$, the interconnection of $G$ and $\tau \Delta$ is well-posed;
		\item for every $\tau \in [0,1]$, the IQC defined by $\Pi$ is satisfied by $\tau \Delta$;
		\item there exists $\epsilon > 0$ such that
		\begin{align}\label{eq:IQC theorem condition 3}
			\begin{bmatrix}
				G (j \omega) \\ I
			\end{bmatrix}^{*}
			\Pi (j\omega)
			\begin{bmatrix}
				G (j \omega) \\ I
			\end{bmatrix}
			\leq - \epsilon I, ~ \forall \omega \in \mathbb{R}.
		\end{align}
	\end{enumerate}
	Then, the interconnection of $G$ and $\Delta$ is stable.
\end{theorem}

Note that if $\Pi (j \omega) = \begin{bmatrix}
	\Pi_{11} (j \omega) & \Pi_{12}(j \omega) \\ \Pi_{12}^* (j \omega) & \Pi_{22} (j \omega)
\end{bmatrix}$ satisfies $\Pi_{11} (j \omega) \geq 0$ and $\Pi_{22} (j \omega) \leq 0$, then the condition $\Delta \in \text{IQC}(\Pi)$ implies that $\tau \Delta \in \text{IQC} (\Pi)$ for all $\tau \in [0, 1]$.

The IQC theorem \il{for discrete-time} systems can be found in, e.g., \cite{jonsson2001lecture}.

\subsection{Mirror descent algorithm}
Consider the optimization problem
\begin{align}\label{eq:optimization problem}
\min_{x \in \mathcal{X}} f(x)
\end{align}
where $\mathcal{X}$ is \icl{a} closed and convex constraint set and $ \mathcal{X} \subseteq \mathbb{R}^{d}$, $f$ is the objective function and $f \in S(\mu, L)$.
For simplicity, We will consider the unconstrained case in this work first, i.e., $\mathcal{X} = \mathbb{R}^{d}$, 
and extend the results to constraint set in the future.
%

We can solve \eqref{eq:optimization problem} with the well-known gradient descent (GD) algorithm
$
x_{k+1} = x_k - \eta \nabla f(x_k),
$
or equivalently,
\begin{align*}
x_{k+1} = \underset{x \in \mathbb{R}^{d}}{\text{argmin}} \left\{ \nabla f(x_k)^T  x  + \frac{1}{2 \eta} \| x - x_{k} \|^2_2  \right\}
\end{align*}
where $\eta > 0$ is a fixed stepsize. Observe that the Euclidean norm used above can be replaced with other \icl{distance \il{measures}}
to generate new algorithms.

The Bregman divergence defined with respect to a distance generating function (DGF) $\phi: \mathbb{R}^{d} \rightarrow \mathbb{R}$ is given by
\begin{align}
D_{\phi} (y, x) = \phi (y) - \phi (x) - ( y - x )^T \nabla \phi (x).
\end{align}
where $\phi(x) \in S (\mu_{\phi}, L_{\phi})$.
Then, the MD algorithm is given by
\begin{align}\label{eq:MD algorithm compact}
x_{k+1} = \underset{x \in \mathbb{R}^{d}}{\text{argmin}} \left\{ \nabla f(x_k)^T x + \frac{1}{\eta} D_{\phi} (x, x_{k}) \right\}.
\end{align}
Denote $\bar{\phi}$ as the convex conjugate of function $\phi$, i.e.,
\begin{align*}
\bar{\phi}(z) = \sup_{x} \left\{ x^T z - \phi (x) \right\}.
\end{align*}
Denote $\mu_{\bar{\phi}} = L_{\phi}^{-1}$, and $L_{\bar{\phi}} = \mu_{\phi}^{-1} $. It follows that $\bar{\phi} \in S (\mu_{\bar{\phi}}, L_{\bar{\phi}} )$, and
$
z = \nabla \phi(x) \Longleftrightarrow x = \nabla \bar{\phi} (z).
$
In other words, $\nabla \bar{\phi}$ is the inverse function of $\nabla \phi$.
Then, the MD algorithm \eqref{eq:MD algorithm compact} can be written as
\begin{align*}
z_{k+1}  = z_{k} - \eta \nabla f (x_{k}), \quad x_{k+1} = \nabla \bar{\phi} (z_{k+1} )
\end{align*}
or equivalently,
\begin{align}\label{eq:MD discrete-time composition}
z_{k+1} = z_{k} - \eta ( \nabla f  \circ \nabla \bar{\phi} )(z_{k})
\end{align}
where $\circ$ represents \il{composition of functions}.
Similarly, the continuous MD algorithm can be given by
\begin{align}\label{eq:MD continuous-time composition}
\dot{z} \nkl{(t)} =  - \eta (\nabla f  \circ \nabla \bar{\phi})(z\nkl{(t)}).
\end{align}
\ml{Any equilibrium point of the above systems satisfies $\nabla f \left( \nabla \bar{\phi} (z^\text{\nkl{opt}}) \right) = \nabla f(x^\text{\nkl{opt}}) = 0_d$, which is the optimal solution to problem \eqref{eq:optimization problem}.
}

\nkl{In the remainder of this paper, the time dependency in the continuous-time case
will be omitted to simplify the notation.}

\todoing{\st{how is $\phi$ chosen? Can it be arbitrary? Why is the equilibrium point unchanged?}}
\ml{Note that the DGF $\phi$ can be an arbitrary function in $S (\mu_{\phi}, L_{\phi})$. 
\il{Function} $\phi$ is \il{usually} chosen such that its convex conjugate is easily computable.
The principal motivation is to generate a distance function that reflects the geometry of the given constraint set $\mathcal{X}$ so that it can often be automatically eliminated during calculation.
Various examples such as minimization over the unit simplex via the Kullback-Leibler divergence can be found in \cite{beck2003mirror,bubeck2014convex,sun2021centralized} and references therein. }

\section{Continuous-time mirror descent method}\label{Continuous-time mirror descent method}
In this section, we construct an IQC framework to analyze the continuous-time MD method.

\subsection{MD algorithm in the form of Lur'e systems}
\ml{It seems that the composition of operators in \eqref{eq:MD continuous-time composition} hinders the \il{direct} application of \il{an} IQC framework since the \icl{composite} operator may not belong to the original classes of the two operators, e.g., the composition of two monotone operators is not necessarily monotone.}
Nevertheless, the cascade connection of two nonlinear operators can be transformed into the feedback interconnection of a linear system with the direct sum of the two nonlinear operators, similarly to the example in \cite{megretski1997system}.
Therefore, the continuous-time MD algorithm \eqref{eq:MD continuous-time composition} can be rewritten as
\begin{align}\label{eq:linear system continuous-time}
\dot{z}  = A z + B u, \quad
y =
C z +
D u
\end{align}
where $u = \begin{bmatrix} u_1\\ u_2 \end{bmatrix}$, $y = \begin{bmatrix}
y_1 \\ y_2
\end{bmatrix}$, the system matrices are
\begin{align}\label{eq:system matrices continuous-time}
\begin{bmatrix}
	\begin{array}{c|c}
      A & B\\
      \hline
      C & D
    \end{array}
\end{bmatrix}
=
\begin{bmatrix}
	\begin{array}{c|cc}
		-\eta \mu_f \mu_{\bar{\phi}} I_d &
		-\eta I_d & -\eta \mu_f I_d \\
		\hline
		\mu_{\bar{\phi}} I_d & 0_d & I_d\\
		I_d & 0_d & 0_d
		\end{array}
\end{bmatrix}
\end{align}
and the system input is
\begin{align}\label{eq:system input continuous-time}
\begin{bmatrix}
u_1 \\ u_2
\end{bmatrix}
= &
\begin{bmatrix}
\nabla f (y_1) - \mu_{f} y_1 \\ \nabla \bar{\phi} (y_2) - \mu_{\bar{\phi}} y_2
\end{bmatrix}.
\end{align}
The transfer function matrix of the linear system is
\begin{equation}\label{eq:transfer function continuous-time}
\begin{aligned}
G(s) = & C(sI_{d} - A)^{-1} B + D\\
= &
\frac{1}{s + \eta \mu_{f} \mu_{\bar{\phi}} }
\begin{bmatrix}
- \eta \mu_{\bar{\phi}} & s\\
- \eta & - \eta \mu_{f}
\end{bmatrix} \otimes I_d
\end{aligned}
\end{equation}
\lmm{where $\otimes$ denotes the Kronecker product.}

Next, define $z^\text{\nkl{opt}}$, $x^\text{\nkl{opt}}$ as the \nkl{optimal  state} with corresponding $x^\text{\nkl{opt}}$, $y^\text{\nkl{opt}}$ and $u^\text{\nkl{opt}}$.
Let $\tilde{z} = z - z^\text{\nkl{opt}}$, $\tilde{y} = y - y^\text{\nkl{opt}}$, $\tilde{u} = u - u^\text{\nkl{opt}}$.
We obtain the error system
\begin{equation}\label{eq:mirror descent algorithm continuous-time}
\dot{\tilde{z}} =   A \tilde{z} + B \tilde{u}, \quad
\tilde{y} =
C
\tilde{z} +
D
\tilde{u}\end{equation} with \begin{equation}
\begin{split}
\tilde{u}
:= & {\Delta} \left( \begin{bmatrix} y_1 - y^\text{\nkl{opt}}_1\\ y_2 - y^\text{\nkl{opt}}_2 \end{bmatrix} \right)
 =  \begin{bmatrix} \Delta_1 \left(y_1 - y^\text{\nkl{opt}}_1 \right) \\ \Delta_2 \left( y_2 - y^\text{\nkl{opt}}_2 \right)\end{bmatrix},
\end{split}
\label{eq:nonlinearity}
\end{equation}
where $\Delta_1 (x)$, $\Delta_2 (x)$ are defined by
$$
\Delta_1 (x) \hspace{-1mm} = \hspace{-1mm} \left(\nabla f (x + y_1^{\textup{opt}}) \hspace{-1mm} - \hspace{-1mm} \mu_{f} (x + y_1^{\textup{opt}}) \right) - \left(\nabla f (y^\text{\nkl{opt}}_1) \hspace{-1mm} - \hspace{-1mm} \mu_{f} y^\text{\nkl{opt}}_1 \right)
$$
$$
\Delta_2 (x) \hspace{-1mm} = \hspace{-1mm} \left(\nabla \bar{\phi} (x + y_2^{\textup{opt}}) \hspace{-1mm} - \hspace{-1mm} \mu_{\bar{\phi}} (x + y_2^{\textup{opt}}) \right) - \left(\nabla \bar{\phi} (y^\text{\nkl{opt}}_2) \hspace{-1mm} - \hspace{-1mm} \mu_{\bar{\phi}} y^\text{\nkl{opt}}_2 \right).
$$

\begin{figure}[htbp]
\centering
\subfigure[]{\includegraphics[width = 0.8\linewidth]{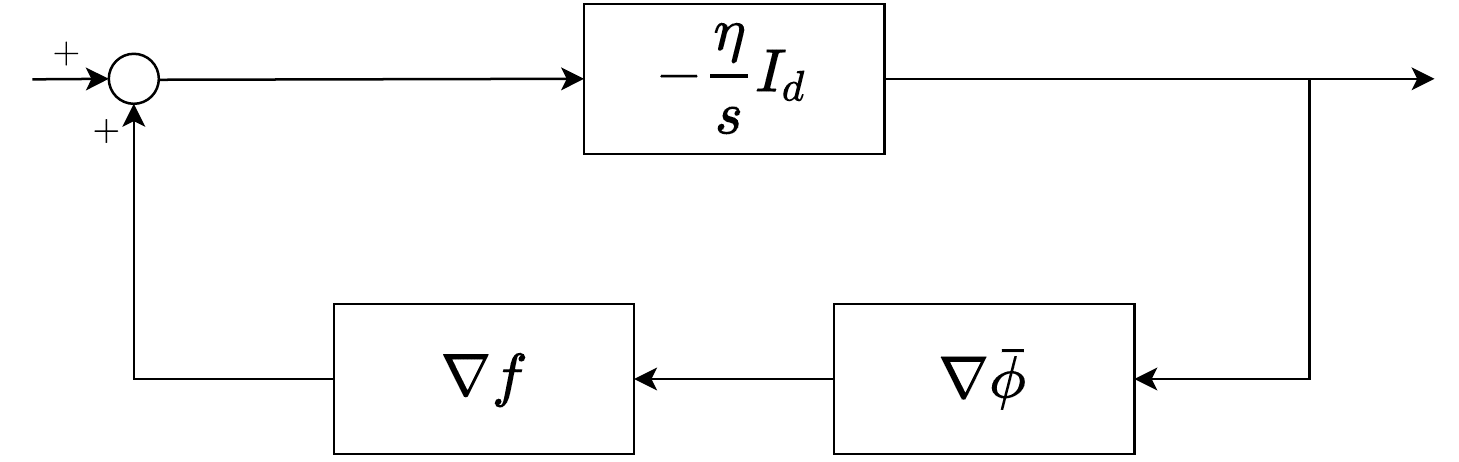}\label{fig: feedback cascade}}
\subfigure[]{\includegraphics[width = 0.8\linewidth]{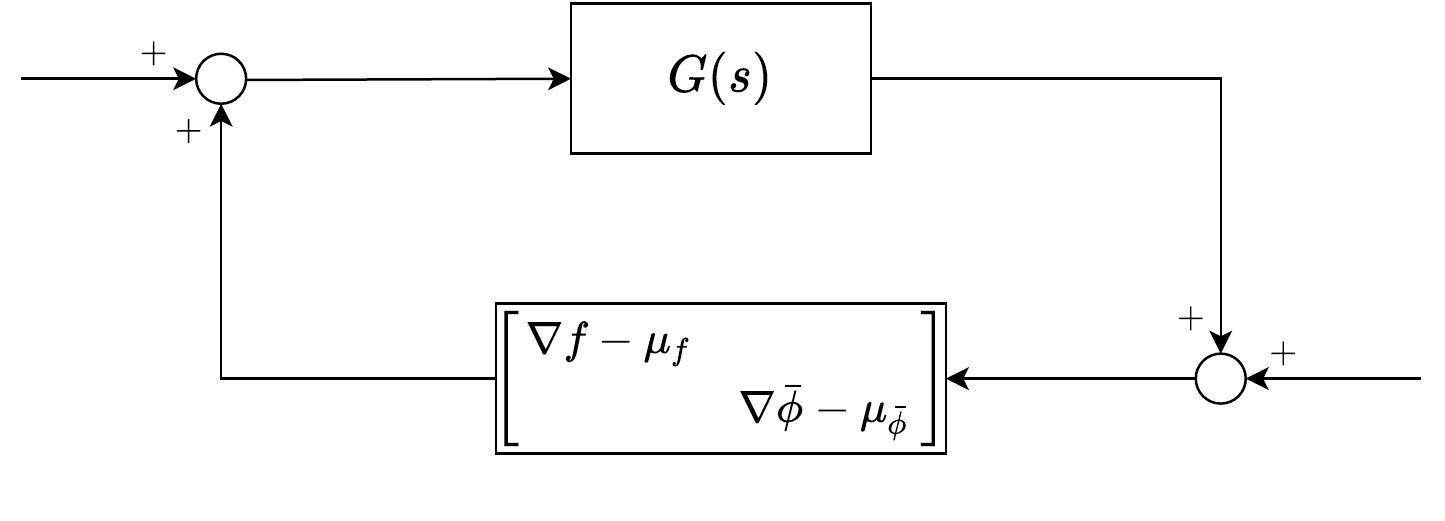}\label{fig: feedback_cascade_equivalence}}
\subfigure[]{\includegraphics[width = 0.8\linewidth]{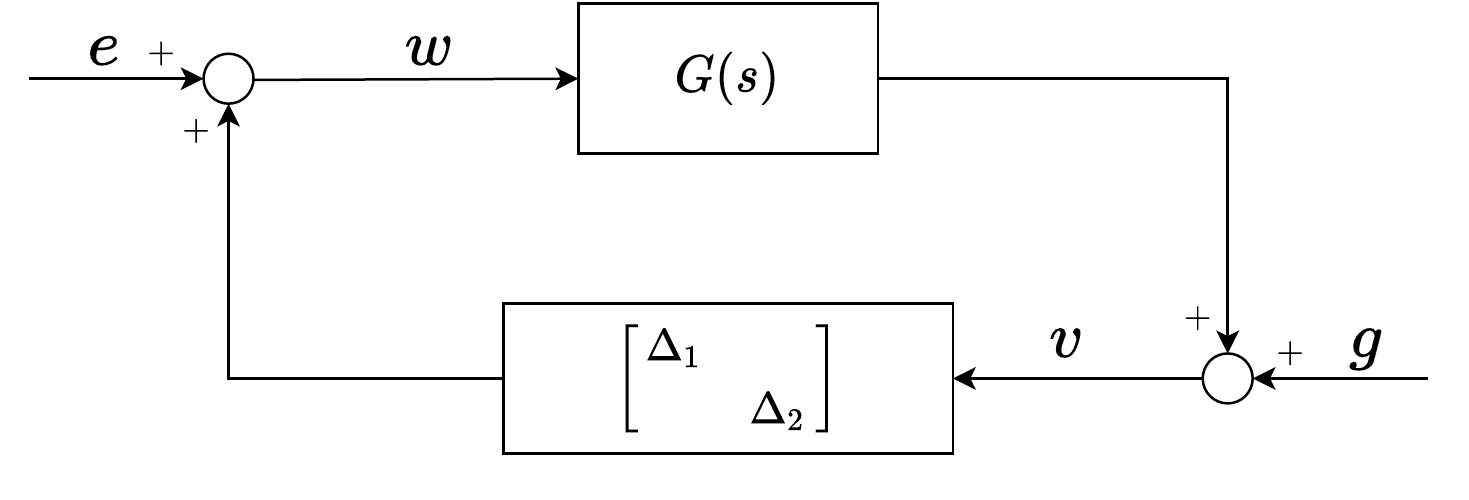}\label{fig: error_system}}
\caption{\ml{Transformation of the MD method to \il{a Lur'e} system. (a) represents the composition of operators, which is transformed to the direct sum of operators in (b), where $G(s)$ is given by \eqref{eq:transfer function continuous-time}.
(c) is the error system in \eqref{eq:mirror descent algorithm continuous-time}, and $\Delta_1$, $\Delta_2$ are given by \eqref{eq:nonlinearity}. }}
\label{fig:feedback interconnection of nonlinear functions}
\end{figure}

\ml{It is apparent that the above error system is in the form of \il{a Lur'e system} \eqref{eq:feedback interconnection model}, where $v = \tilde{y}$, $w = \tilde{u}$, $e = 0_d$, and $g$ \il{is a trajectory that} represents the \il{effect of the} initial condition. The transformation can be depicted by Fig.~\ref{fig:feedback interconnection of nonlinear functions}.}

\subsection{IQCs for gradients of convex functions}
\todoiny{\st{I think more text is needed before starting section III.B to justify why we are bringing all these materials.}}
\ml{In this subsection, we will include a group of useful IQCs for gradients of convex functions to characterize the nonlinearity $\Delta$.
Note that conic combinations of various IQCs are also valid IQCs which better characterize the nonlinearity and lead to less conservative stability margins.
\todoiny{\st{Say something about the conic combination}}
}
\subsubsection{Sector IQC}
\ml{The sector IQC is introduced in the following lemma as a result of the co-coercivity of gradients.}
\begin{lemma}[\hspace{1sp}\cite{lessard2016analysis}]\label{lem:co-coercivity sector IQC}
Suppose \lmm{a function $f \in S(\mu, L)$}. For all $x, y$, the following quadratic constraint (QC) is satisfied,
\begin{align*}
\left[
\begin{smallmatrix}
y - x \\
\nabla \lmm{f} (y) - \nabla  \lmm{f} (x)
\end{smallmatrix}
\right]^T
\left[
\begin{smallmatrix}
-2 \mu L I_d & (L + \mu) I_d \\
(L + \mu )I_d & -2 I_d
\end{smallmatrix}
\right]
\left[
\begin{smallmatrix}
y - x \\
\nabla \lmm{f} (y) - \nabla \lmm{f} (x)
\end{smallmatrix}
\right]
\geq 0.
\end{align*}
\end{lemma}

\nkl{Note that as  $f \in S(\mu_f, L_f)$, $\bar{\phi} \in S(\mu_{\bar{\phi}}, L_{\bar{\phi}})$, \ml{then $f (\cdot) - \frac{1}{2} \mu_f \| \cdot \|^2 \in S ( 0, L_f - \mu_f )$ and $\bar{\phi} (\cdot) - \frac{1}{2} \mu_{\bar{\phi}} \|\cdot \|^2 \in S ( 0, L_{\bar{\phi}}- \mu_{\bar{\phi}} )$.}
\todoinc{are you meant to substract just a constant to the functions? what is the role of $\|y_i\|^2$, and the division by 2?}
Moreover, using~\Cref{lem:co-coercivity sector IQC},   $\Delta \in \text{IQC}( \Pi_{s})$, where $\Delta$ is defined in \ml{\eqref{eq:nonlinearity}}  and}
\begin{align}\label{eq:IQC sector bounded}
\Pi_{s} \hspace{-1mm} =
\left[
\begin{smallmatrix}
0_d & 0_d & \alpha_1 ( L_f - \mu_f )  I_d & 0_d\\
0_d & 0_d & 0_d & \alpha_2 ( L_{\bar{\phi}} - \mu_{\bar{\phi}} ) I_d\\
\alpha_1 ( L_f - \mu_f )  I_d & 0_d & -2 \alpha_1 I_d & 0_d\\
0_d & \alpha_2 ( L_{\bar{\phi}} - \mu_{\bar{\phi}} ) I_d & 0_d & -2 \alpha_2 I_d
\end{smallmatrix}
\right]
\end{align}
\nkl{where 
 $\alpha_1, \alpha_2 \geq0$.}

\subsubsection{Popov IQC}
\ml{The Popov IQC is introduced as follows.}
\begin{lemma}\label{lem:popov IQC}
Suppose $f \in S(0, L)$. The nonlinearity $\nabla f(x) - \nabla f(x^\textup{\nkl{opt}})$ satisfies the Popov IQC by $\Pi_{P} (j \omega)$ given by
\begin{align*}
\Pi_{P} (j \omega) =
	\pm \begin{bmatrix}
			0_d & - j \omega I _d \\ j\omega I_d & 0_d
		\end{bmatrix}.
\end{align*}
\end{lemma}
As $f (\cdot) - \frac{1}{2} \mu_f \| \cdot \|^2 \in S ( 0, L_f - \mu_f )$ and $\bar{\phi} (\cdot) - \frac{1}{2} \mu_{\bar{\phi}} \| \cdot \|^2 \in S ( 0, L_{\bar{\phi}}- \mu_{\bar{\phi}} )$,
using \Cref{lem:popov IQC}, we have $\Delta \in \text{IQC} (\Pi_{p} \lmm{(j\omega)})$, where $\Delta$ is defined in \eqref{eq:nonlinearity} and 
\begin{align}\label{eq:Popov IQC for delta}
\Pi_{p}(j \omega) = 
\left[
\begin{smallmatrix}
0_d & 0_d & -j \omega \beta_1  I_d & 0_d\\
0_d & 0_d & 0_d & - j \omega \beta_2 I_d\\
 j \omega \beta_1  I_d & 0_d & 0_d & 0_d\\
0_d & j \omega \beta_2 I_d & 0_d & 0_d
\end{smallmatrix}
\right]
\end{align}
where $\beta_1, \beta_2 \geq 0$.


%

\subsection{Convergence analysis via \il{IQCs} in \nkl{frequency domain}}
\todoing{\st{Need some introductory text to summarize the convergence results that will be presented.}}

In this subsection, we will present the convergence analysis of the MD method.
There is a rich literature showing the convergence of the MD method, e.g.,\cite{nemirovskij1983problem,beck2003mirror,krichene2015accelerated}. We show that using \il{an} IQC analysis also leads to such a conclusion.

\begin{theorem}\label{thm:convergence continuous-time IQC}
	Consider \icl{a} Lur'e system \icl{described by} \eqref{eq:feedback interconnection model} where $g \in \mathbf{L}_{2}[0, \infty)$, $e \in \mathbf{L}_{2}[0, \infty)$, $G(s)$ is given by \eqref{eq:transfer function continuous-time}, \ml{$\Delta$ is defined in \eqref{eq:nonlinearity}} with $f \in S(\mu_f, L_f)$, $\phi \in S(\mu_\phi, L_\phi)$.
The system is stable and \ml{the trajectory of $x = \nabla \bar{\phi} (z)$ with any 
\il{initial} condition $z(0) = z_0$ of the MD method \eqref{eq:MD continuous-time composition} converges to the optimal solution of problem \eqref{eq:optimization problem}.}
\todoing{\st{Need to be careful about the theorem statement - stability in IQCs is in an input/output sense, whereas here you want to deduce convergence for any initial conditions; see e.g. how Lessard paper handles this.}}
\end{theorem}
\lk{\textbf{Sketch of the proof.} 
Stability  can be shown using
 \Cref{thm:IQC theorem} with $\Pi ( j \omega )$ given~by \begin{align}\label{eq:IQC sector bounded + Popov}
		\Pi ( j \omega ) = &
		\left[
		\begin{smallmatrix}
			0_d & 0_d & \left( \alpha_1 ( L_f - \mu_f ) - \beta_1 j \omega \right) I_d & 0_d \\
			0_d & 0_d & 0_d & \left( \alpha_2 ( L_{\bar{\phi}} - \mu_{\bar{\phi}} ) - \beta_2 j \omega \right) I_d\\
			* & * & -2 \alpha_1 I_d & 0_d \\
			* & * & 0_d & -2 \alpha_2 I_d
		\end{smallmatrix}
		\right].
	\end{align}
	Note that $\Pi(j\omega)$ of~\eqref{eq:IQC sector bounded + Popov} is the conic combination of the sector IQC \eqref{eq:IQC sector bounded} and Popov IQC \eqref{eq:Popov IQC for delta}. The stability implies that $v \rightarrow 0$ as $t \rightarrow \infty$,  
	which means the trajectories of the error system~\eqref{eq:mirror descent algorithm continuous-time} tend to $0$ as $t \rightarrow \infty$ and  the trajectory of $x = \nabla \bar{\phi} (z)$ for {any input} $z_0$ converges to the optimal solution of problem \eqref{eq:optimization problem}. \qed } 

\Cref{thm:convergence continuous-time IQC} is based on conditions in the frequency domain, which \icl{do not} describe the convergence rate of the MD algorithm. To this end, we will \icl{investigate, in the next subsection,} the MD method in the time domain and reveal the connection between the Bregman divergence function and the Popov \icl{criterion.}

\subsection{\nkl{Convergence analysis via \il{IQCs} in time domain}}

In this subsection, we show that the Bregman divergence function, \icl{which is widely used as a Lyapunov function for the MD algorithm, is a special case of Lyapunov functions that are associated with the Popov criterion.} \todoing{this sentence is vague}
\icl{This connection is established by applying the multivariable Popov criterion,
which is \il{adapted} 
from} \cite{moore1968generalization,khalil1996nonlinear,jonsson1997stability,carrasco2013equivalence}.
\begin{lemma}\label{lem:popov}
\nkl{Let $H (s) \in \mathbf{RH}_{\infty}^{p\times p}$ and let $\psi  : \mathbb{R}^{p} \to \mathbb{R}^{p}$ be a memoryless nonlinearity composed of $p$ memoryless nonlinearities $\psi_i$
 with each   being }\ml{slope-restricted on sector \textup{[0, $k_i$]}}, i.e., \ml{$0 \leq \frac{\psi_i(x_1) - \psi_i(x_2)}{x_1 - x_2} \leq k_i$, $\forall x_1 \neq x_2$}, $0< k_i \leq \infty$, for $1\leq i \leq p$. If there exist constants $\alpha_i \geq 0$ and $\gamma_i \geq 0$ such that $ \text{Re} \left\{ \alpha K^{-1} + (\alpha + j \omega \Gamma) H (j\omega ) \right\} \geq \delta$ for some $\delta > 0$, where $K = \textup{diag}(k_1, \ldots, k_p)$, $\alpha = \textup{diag} (\alpha_1, \ldots, \alpha_p )$, $\Gamma = \textup{diag}(\gamma_1, \ldots, \gamma_p)$.
Then, the negative feedback interconnection of $H(s)$ and $\psi$ is stable.
\end{lemma}
\todoing{\st{agree with above - what does $\mathbf{L}_2$ stable mean?} }
\begin{remark}
The parameters $\alpha_i$, $\gamma_i$ result directly from the conic parameterization of the sector and Popov IQCs. \ml{The proof of \Cref{thm:convergence continuous-time IQC} can be seen as an application of \Cref{lem:popov} since the third condition in \Cref{thm:IQC theorem}, \il{with the IQC used in the proof of \Cref{thm:convergence continuous-time IQC},}  is equivalent to the inequality condition \il{in \Cref{lem:popov}.}} 
It is noteworthy that the consideration of $\alpha_i$ is crucial since it provides \kl{more flexibility and thus less conservatice} results for the MIMO \il{case \cite{moore1968generalization}.} 
The original Popov criterion requires that the linear system $H(s)$ is strictly proper, i.e., there is no direct feedthrough term \cite{moore1968generalization,khalil1996nonlinear}, and the derivative of the input to $\psi$ is bounded \cite{jonsson1997stability}. These restrictions are removed in \cite{carrasco2013equivalence}.

\todoing{\st{as above $\mathbf{L}_2$ stable needs clarification}}
\end{remark}

\todoiny{\st{More details are needed.}}

We can apply \Cref{lem:popov} and obtain a condition to characterize the exponential convergence rate for the continuous-time MD method.
\begin{theorem}\label{thm:convergence rate continuous-time}
The continuous-time MD algorithm \eqref{eq:MD continuous-time composition} with $f \in S(\mu_f, L_f)$, $\phi \in S(\mu_\phi, L_\phi)$, converges exponentially to the optimal solution \nkl{with a convergence rate $\rho$} 
if there exist $P > 0$, $\Gamma = \textup{diag}(0, \gamma) \geq 0$, and $\alpha = \textup{diag} ( \alpha_1, \alpha_2) \geq 0$ such that
\begin{align}\label{eq:LMI exponential rate continuous-time}
	\begin{bmatrix}
		P\tilde{A} + \tilde{A}^TP + 2 \rho P & P \tilde{B} - \tilde{C}^T \\
		* & - \left( \tilde{D} + \tilde{D}^T \right)
	\end{bmatrix}
	\leq 0
\end{align}
where $\tilde{A} = A$, $\tilde{B} = - B$, $\tilde{C} = \lmm{(\alpha + \rho \Gamma)}  C + \Gamma C A$, $\tilde{D} = - \alpha D + \alpha K^{-1} - \Gamma C B$ and $(A,B,C,D)$ is defined in \eqref{eq:system matrices continuous-time}.
\end{theorem}
\lmm{\textbf{Sketch of the proof.} Apply \Cref{lem:popov} with $H(s) = - G(s)$ and $\psi = \Delta$, and the Kalman-Yakubovich-Popov (KYP) Lemma\cite{rantzer1996kalman}, taking into account the exponential stability\cite[Theorem 2]{hu2016exponential}. \qed
}

The convergence rate $\rho$ in \eqref{eq:LMI exponential rate continuous-time} needs to be treated as a \icl{constant} 
such that \eqref{eq:LMI exponential rate continuous-time} is an LMI. Nevertheless, a bisection search on $\rho$  can be carried out in \eqref{eq:LMI exponential rate continuous-time} to obtain the largest admissible convergence rate for the continuous-time MD \eqref{eq:MD continuous-time composition}.

\begin{remark}
\lmm{\Cref{thm:convergence rate continuous-time}  follows from an application of the \lk{multivariable}  Popov criterion and its corresponding Lyapunov function which is
\begin{equation}\label{eq:lyapunov function popov}
\begin{aligned}
V = & \frac{1}{2} \tilde{z}^T P \tilde{z} + \gamma \int_{y^\text{\nkl{opt}}_2}^{y_2} \psi_2 (\tau)  d \tau\\
=  & \frac{1}{2} \tilde{z}^T P\tilde{z} + \gamma D_{\bar{\phi}}(z, z^\text{\nkl{opt}}) - \frac{\gamma \mu_{\bar{\phi}}}{2} \| \tilde{z} \|^2.
\end{aligned}
\end{equation}}
When $P = \lmm{\gamma \mu_{\bar{\phi}} I_d}$ and $\gamma = 1$, the Lyapunov function \eqref{eq:lyapunov function popov} reduces to the Bregman divergence function, which is 
\il{a common} 
choice of Lyapunov function for the MD method \cite{nemirovskij1983problem, krichene2015accelerated}. This implies that in the analysis of convergence rate, using the IQC analysis framework with a conic combination of \ml{IQCs} including the Popov and Zames-Falb-O'Shea \todoing{\st{this is a special class of Zames-Falb multiplier - right?}} \il{ones, yields} an equivalent or less conservative worst-case convergence rate, \il{to the one that follows by simply} using the Bregman-type Lyapunov functions.
\end{remark}
\todoing{A question that might come up is why we are using only the Popov criterion in Theorem \ref{thm:convergence rate continuous-time} and not more general Zames Falb multipliers. Also Zames-Falb multipliers can be non-causal in which case a Lyapunov function cannot be used and the convergence rate proof needs to be done in a different way, as in e.g. [22].\\
Answer:
It is true that the convergence proof may be different if Zames-Falb multipliers are used.
However, in this work, using the Popov IQC is sufficient to obtain a tight convergence rate. That is, there is no need to apply Zames-Falb multipliers for this case.
}

\section{Discrete-time mirror descent method}\label{Discrete-time mirror descent method}
Similar to the continuous-time case, the discrete-time MD \il{algorithm} in \eqref{eq:MD discrete-time composition} can be rewritten into the following Lur'e system,
\begin{align}\label{eq:mirror descent discrete-time}
	z_{k+1} = A z_{k} + B u_k, \quad y_k = C z_k + D u_k
\end{align}
where $ u_k = \begin{bmatrix} u_k^{(1)} \\ u_k^{(2)} \end{bmatrix}$, $y_k = \begin{bmatrix} y_k^{(1)}\\ y_k^{(2)} \end{bmatrix}$ the system matrices are
\begin{align}\label{eq:system matrices discrete-time}
\hspace{-2mm}
\begin{bmatrix}
	\begin{array}{c|c}
      A & B\\
      \hline
      C & D
    \end{array}
\end{bmatrix}
=
\begin{bmatrix}
	\begin{array}{c|cc}
		(1-\eta \mu_f \mu_{\bar{\phi}})I_d &
		-\eta I_d & -\eta \mu_f I_d \\
		\hline
		\mu_{\bar{\phi}} I_d & 0_d & I_d\\
		I_d & 0_d & 0_d
		\end{array}
\end{bmatrix}
\end{align}
and the system input is
\begin{align}\label{eq:input discrete-time}
	\begin{bmatrix}
u_k^{(1)} \\ u_k^{(2)}
\end{bmatrix}
= &
\begin{bmatrix}
\nabla f (y_k^{(1)}) - \mu_{f} y_k^{(1)} \\ \nabla \bar{\phi} (y_k^{(2)}) - \mu_{\bar{\phi}} y_k^{(2)}
\end{bmatrix}
.
\end{align}
\ml{Defined $z^\textup{\nkl{opt}}$ as the optimal \il{value of $z$ at} steady state, with corresponding \il{equilibrium values} $y^\textup{\nkl{opt}} = \left( y^{(1),\nkl{\text{opt}}}, y^{(2),\nkl{\textup{opt}}}\right)$, $x^\textup{\nkl{opt}}$, and $u^\textup{\nkl{opt}}$. Define $\tilde{u}_{k} = u_k - u^\textup{\nkl{opt}}$, then we have}
\begin{align}\label{eq:nonlinearity discrete-time}
\tilde{u}_{k}
= \Delta
\left(
\begin{bmatrix}
	y_k^{(1)} - y^{(1),\nkl{\textup{opt}}} \\ y_k^{(2)} - y^{(2),\nkl{\textup{opt}}}
\end{bmatrix}
\right)
\end{align}
\ml{
where the nonlinear operator $\Delta$ in \eqref{eq:nonlinearity discrete-time} is \icl{the same as that used for the continuous-time algorithm in \eqref{eq:nonlinearity}.}}
%

\subsection{Convergence rate via IQC}
There is no exact counterpart for the Popov criterion in \icl{discrete-time}. Similar ones are the Jury-Lee criteria \cite{jury1964stability,haddad1994parameter}.
Though we could easily provide an LMI condition for the discrete-time system \eqref{eq:system matrices discrete-time}, \eqref{eq:nonlinearity discrete-time} following the discrete-time Jury-Lee criteria via the same Lyapunov function, \ml{we remark that in discrete-time, all IQCs to characterize monotone and bounded nonlinearities are within the set of Zames-Falb-O'Shea IQCs.}
\todoing{\st{do not fully understand the previous sentence. What do you mean by saying the multipliers "preserve the positivity of monotone and bounded nonlinearities".}}
Therefore, we can directly apply the class of 
 Zames-Falb-O'Shea \il{IQCs with a state-space representation as in} \cite{lessard2016analysis}. We will only adopt a simple type of the Zames-Falb-O'Shea IQC here \ml{because \il{this} is sufficient to obtain a tight convergence rate for the MD method.}
\todoing{\st{is this done for simplicity in the presentation?}\\
Same as the last answer. The simplest type is sufficient to obtain a tight convergence rate.
}
\todoKL{$\Pi_{w}$ does not depend on $j \omega$.}

From \cite{lessard2016analysis}, we can obtain that $\Delta$ satisfies the weighted-off-by-one IQC defined by
$
		\Pi_{w} =  \Psi_{w}^{*} M_{w} \Psi_{w}, \quad M_{w} = \begin{bmatrix} 0_{2d} & \beta I_{2d} \\ \beta I_{2d} & 0_{2d}\end{bmatrix},
$
\ml{where $\Psi_{w}$ is a transfer function matrix with the following state-space representation,}
	\begin{align}
	\hspace{-2mm}
	\begin{bmatrix}
		\hspace{-0.5mm}
			\begin{array}{c|cc}
				A_{\Psi_{w}} & B_{\Psi_{w}}^{y} & B_{\Psi_{w}}^{u} \\
				\hline
				C_{\Psi_{w}} & D_{\Psi_{w}}^{y} & D_{\Psi_{w}}^{u}
			\end{array}
		\hspace{-0.5mm}
		\end{bmatrix} \hspace{-0.5mm} = \hspace{-0.5mm}
		\begin{bmatrix}
		\hspace{-0.5mm}
			\begin{array}{c|cc}
				0_{2d} & - K I_{2d} & I_{2d} \\
				\hline
				\bar{\rho}^2 I_{2d} & K I_{2d} & -I_{2d}\\
				0_{2d} & 0_{2d} & I_{2d}
			\end{array}
		\hspace{-0.5mm}
		\end{bmatrix} \hspace{-1mm}
	\end{align}
with $K = \textup{diag}\{L_f - \mu_f, L_{\bar{\phi}} - \mu_{\bar{\phi}}\} \otimes I_{d}$, $\beta = \textup{diag} \{\beta_1, \beta_2 \} \otimes I_{d} \geq 0$, and $\bar{\rho} \geq 0$.

From \Cref{lem:co-coercivity sector IQC}, we can obtain that $\Delta$ satisfies the IQC defined by
$
 \Pi_{s} = \Psi_{s}^* M_{s} \Psi_s, \quad M_s = \begin{bmatrix} 0_{2d} & \alpha I_{2d} \\ \alpha I_{2d} & 0_{2d}\end{bmatrix},
$
\ml{where $\Psi_{s}$ is a transfer function matrix with the following state-space representation,}
\begin{align}
\begin{bmatrix}
	\begin{array}{c|cc}
		A_{\Psi_{s}} & B_{\Psi_{s}}^{y} & B_{\Psi_{s}}^{u}\\
		\hline
		C_{\Psi_{s}} & D_{\Psi_{s}}^{y} & D_{\Psi_{s}}^{u}
		\end{array}
\end{bmatrix} = 
	\begin{bmatrix}
	\begin{array}{c|cc}
		0_{2d} & 0_{2d} & 0_{2d} \\
		\hline
		0_{2d} & K I_{2d} & -I_{2d}\\
		0_{2d} & 0_{2d} & I_{2d}
		\end{array}
\end{bmatrix}
\end{align}
with $\alpha = \textup{diag} \{ \alpha_1, \alpha_2 \} \otimes I_{d} \geq 0$.

Then, we can characterize the convergence rate for the discrete-time MD method by applying the discrete-time IQC theorem.

\begin{theorem}\label{thm:convergence rate discrete-time}
	The discrete-time MD algorithm \eqref{eq:MD discrete-time composition} with $f \in S(\mu_f, L_f)$ and $\phi \in S(\mu_\phi, L_\phi)$ converges \il{with} a rate $\bar{\rho} \leq \rho \leq 1$ if the following LMI is feasible for some $P > 0$, $\alpha \geq 0$, and $\beta \geq 0$ \nkl{such that}
\begin{align}\label{eq:LMI IQC discrete-time}
	\left[
	\begin{smallmatrix}
		\hat{A}^T P \hat{A} - \rho^2 P & \hat{A}^T P \hat{B}\\
		* & \hat{B}^T P \hat{B}
	\end{smallmatrix} \right]
	+
	\left[\begin{smallmatrix}
		\hat{C} & \hat{D}
	\end{smallmatrix} \right]^T
	\left[\begin{smallmatrix}
		M_{s} & \\ & M_{w}
	\end{smallmatrix}\right]
	\left[
	\begin{smallmatrix}
		\hat{C} & \hat{D}
	\end{smallmatrix} \right]
	\leq 0
	\end{align}
where
\begin{align}
	\hat{A} =
	\left[\begin{smallmatrix}
		A & 0_{d \times 2d} & 0_{d \times 2d}\\
		B_{{\Psi_{s}}}^{y} C & A_{\Psi_{s}} & 0_{2d \times 2d}\\
		B_{{\Psi_{w}}}^{y} C & 0_{2d \times 2d} &  A_{\Psi_{w}}
	\end{smallmatrix} \right],~
	\hat{B} =
	\left[\begin{smallmatrix}
    B\\ B_{\Psi_{s}}^{y} D + B_{\Psi_{s}}^{u} \\ B_{\Psi_{w}}^{y} D + B_{\Psi_{w}}^{u}
	\end{smallmatrix} \right], \nonumber \\
	\hat{C} =
	\left[\begin{smallmatrix}
		D_{\Psi_{s}}^{y} C & C_{\Psi_{s}} & 0_{4d \times 2d}\\ D_{\Psi_{w}}^{y} C & 0_{4d \times 2d} & C_{\Psi_{w}}
	\end{smallmatrix} \right],~
	\hat{D} =
	\left[\begin{smallmatrix}
		D_{\Psi_{s}}^{y} D +D_{\Psi_{s}}^{u} \\ D_{\Psi_{w}}^{y} D+ D_{\Psi_{w}}^{u}
	\end{smallmatrix} \right].
\end{align}
\end{theorem}
The proof is similar to \cite[Theorem 4]{lessard2016analysis} and is omitted here.

\subsection{Stepsize selection}\label{subsection Stepsize selection}
It is well-known that the optimal fixed stepsize for the GD method $x_{k+1} = x_{k} - \eta \nabla f (x_{k})$ is $\eta = \frac{2}{L_f + \mu_f}$, rendering the \ml{\lmm{smallest upper bound} \il{for the convergence rate}} $\rho = \frac{\kappa_f - 1}{\kappa_f + 1}$ \nkl{\ml{where} $\kappa_f = L_f/\mu_f$.}
\todoinc{is"with" meant to be "when"?\\
$Lf$ is always bigger than $\mu_f$, so we use ``where'' instead.}
\todoing{\st{in what sense is this stepsize optimal? Is this over all nonlinearities that have this Lipschitz/strong convexity constant?}
Answer: The stepsize is optimal in that it provides the lowest upper bound of convergence rates.
Yes, it is over all nonlinearities that have Lipschitz/strong convexity properties.
}
Notice that the MD method has a similar structure to the GD method by changing the gradient into the composition of two \il{functions.} 
Thus, we let the stepsize be $\eta = \frac{2}{L_f L_{\bar{\phi}} + \mu_f \mu_{\bar{\phi}}}$ \ml{\il{which is analogous to the} optimal stepsize for the GD method.}
\todoing{\st{Why? Does this have a form of optimality?}\\
I think it is hard to explain why since here it is simply an analogue to the optimal stepsize for the GD method. Yes, the stepsize is optimal in the sense that it provides the lowest upper bound of convergence rates.}
\il{We will} show numerically in \Cref{Numerical Examples} that the LMI in \eqref{eq:LMI IQC discrete-time} is feasible for $\rho = \frac{\kappa - 1}{\kappa + 1}$, where $\kappa = \kappa_f \kappa_{\bar{\phi}}$, and $\kappa_f$, $\kappa_{\bar{\phi}}$ are the condition numbers of $f$, $\bar{\phi}$, respectively.

\section{Numerical Examples}\label{Numerical Examples}
In this section, we present two numerical examples to illustrate the IQC analysis for the MD method in continuous-time and discrete-time, respectively.
\subsection{Continuous-time MD method}
\begin{figure}[htbp]
	\centering
		\includegraphics[width = 0.8\linewidth]{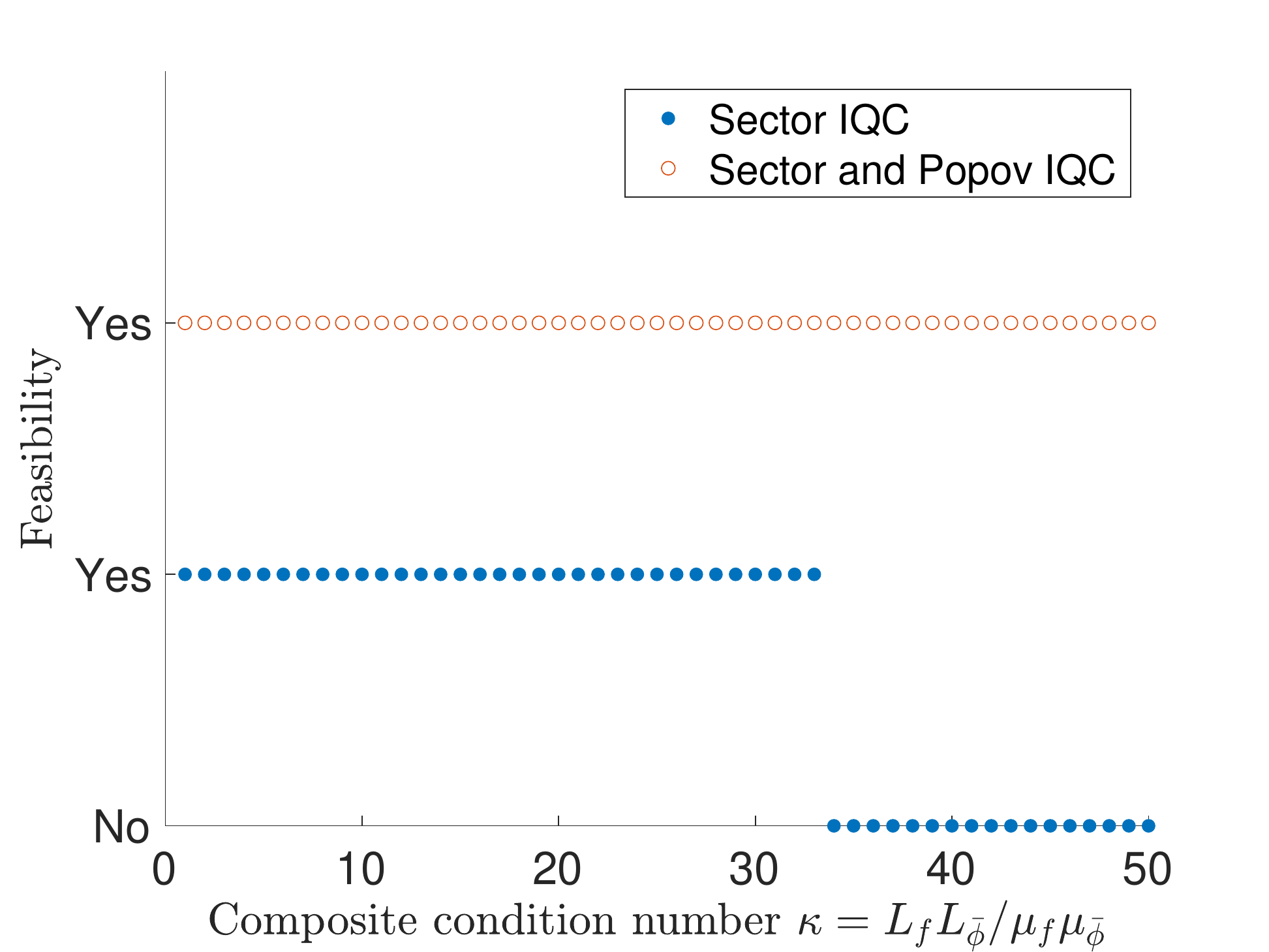}
	\caption{\il{Feasibility of the} problem using \eqref{eq:IQC sector bounded}, the conic combination of  \eqref{eq:IQC sector bounded} and \eqref{eq:Popov IQC for delta}, \il{(for some $\rho>0$) with varying} ratio $\kappa = \frac{L_f L_{\bar{\phi}}}{\mu_{f} \mu_{\bar{\phi}}}$.}
	\label{fig:continuous-time feasibility}
\end{figure}

\todoing{\st{what do you mean by feasibility ranges?}}
\ml{We investigate and compare the \il{feasibility}} of the IQC condition \eqref{eq:IQC theorem condition 3} \il{when} using merely the sector IQC defined by \eqref{eq:IQC sector bounded} and using the conic combination of the sector and Popov IQCs \eqref{eq:IQC sector bounded + Popov}.
\ml{The frequency-\icl{domain} condition \eqref{eq:IQC theorem condition 3} under \eqref{eq:IQC sector bounded} can be easily transformed into \kl{a time-domain condition} via the KYP lemma. While condition \eqref{eq:IQC theorem condition 3} under \eqref{eq:IQC sector bounded + Popov} is satisfied if and only if \eqref{eq:LMI exponential rate continuous-time} in \Cref{thm:convergence rate continuous-time} is feasible for some $\rho > 0$.}
Let $\eta = 1$, $\mu_f = 1$ and $\mu_{\bar{\phi}} = 1$, and $L_f = L_{\bar{\phi}}$.
\il{The feasibility of the IQCs (for some $\rho>0$) with varying composite condition number $\kappa = \frac{L_f L_{\bar{\phi}}}{\mu_{f} \mu_{\bar{\phi}}}$ is} shown in \cref{fig:continuous-time feasibility}.
Note that the MD method should converge for any $L_f > \mu_f$ and $L_{\bar{\phi}} > \mu_{\bar{\phi}}$. However, we can observe that the sector IQC defined by \ml{\eqref{eq:IQC sector bounded}} fails to certify the convergence of the MD method for $\kappa \geq 34$. On the other hand, using the conic combination of the sector IQC \eqref{eq:IQC sector bounded} and the Popov IQC \eqref{eq:Popov IQC for delta}, suffices to certify its convergence for arbitrary $\kappa$.

\subsection{Discrete-time MD method}
\begin{figure}[htbp]
	\centering
	\includegraphics[width = 0.8\linewidth]{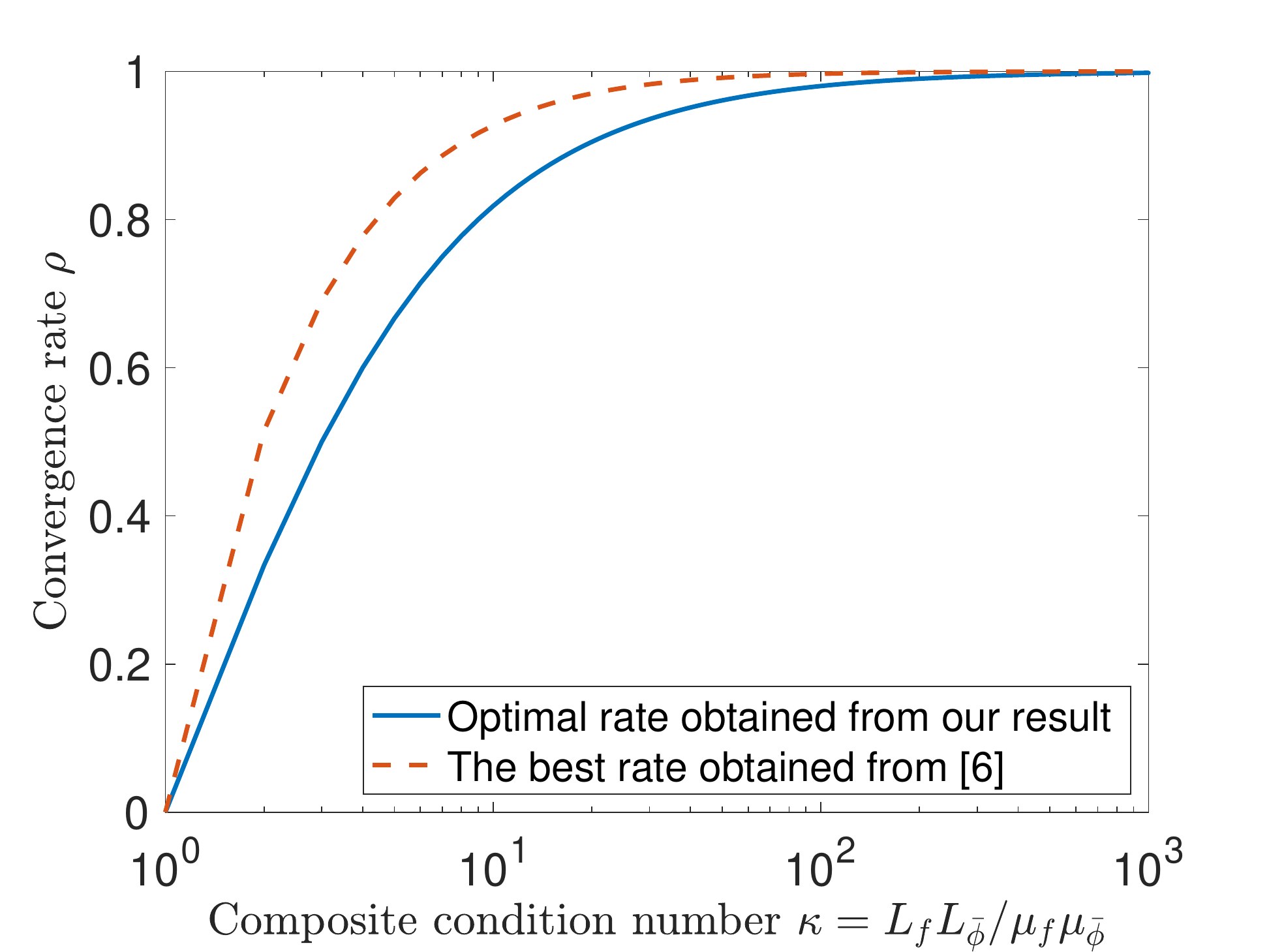}
	\caption{Convergence rate obtained from \eqref{eq:LMI IQC discrete-time} in \Cref{thm:convergence rate discrete-time} and from the SDPs in \cite{sun2021centralized}.  The optimal rate obtained from our result coincides with the curve $\rho =\frac{ \kappa - 1}{\kappa  + 1 }$.}
	\label{fig:discrete-time convergence rate}
\end{figure}

Next, we show the convergence rate for the discrete-time MD method.
Let $\mu_f = 1$, $\mu_{\bar{\phi}} = 1$, and $L_f = L_{\bar{\phi}}$. Let the stepsize be $\eta = \frac{2}{L_f L_{\bar{\phi}} + \mu_f \mu_{\bar{\phi}}}$ as \Cref{subsection Stepsize selection} suggested. We compare the optimal convergence rate obtained from \eqref{eq:LMI IQC discrete-time} in \Cref{thm:convergence rate discrete-time} with that obtained from the SDPs in \cite{sun2021centralized}, where the stepsize and convergence rate are both decision variables. The SDPs in \cite{sun2021centralized} are derived from the Lyapunov function $V (z_k) = \rho^{-k} D_{\bar{\phi}}(z_{k}, z^\textup{\nkl{opt}})$, \todoinc{is $\rho$ in the definition of $V$ the convergence rate - if not use a different symbol?\\ Yes it is the convergence rate.} which \il{is the} Bregman divergence function when $\rho = 1$.
The relation between the composite condition number $\kappa = \frac{L_f L_{\bar{\phi}}}{\mu_{f} \mu_{\bar{\phi}}}$ and the convergence rate $\rho$ is shown in \cref{fig:discrete-time convergence rate}.
We can observe that using the IQC analysis provides a tighter bound \il{for the} convergence rate.
We remark that the convergence rate $\rho =\frac{ \kappa - 1}{\kappa  + 1 }$ obtained here is tight since it is also \ml{the \lmm{smallest upper bound} \il{for the} convergence rates \lmm{of linear systems generated by} all quadratic functions $f \in S(\mu_f, L_f)$ and $\bar{\phi} \in S(\mu_{\bar{\phi}}, L_{\bar{\phi}})$}.
\todoing{\st{why is this the case?}\\
Answer: This can be shown by analyzing the convergence rate of linear systems.}

\section{Conclusion}\label{Conclusion}
An IQC analysis framework has been \il{developed} 
 for the MD method in both continuous-time and \ill{discrete-time. In} \il{continuous-time}, we have shown that the Bregman divergence function is a special case of the \icl{Lyapunov functions associated with the Popov criterion \ill{when these are applied to an appropriate reformulation of the problem}.} In \il{discrete-time,} \ml{we have provided \ill{upper bounds} for the convergence \ill{rate via appropriate IQCs applied to the transformed system. It has also been illustrated via numerical examples that these bounds can be tight.}}
 Future work \il{includes extending the framework developed to other related algorithms such as accelerated MD methods.}

\bibliographystyle{IEEEtran}
\bibliography{References}
\end{document}